\documentclass[12pt,oneside]{amsart}
\usepackage{amsfonts}
\usepackage{amscd}
\usepackage{amsthm}
\usepackage{amssymb}
\usepackage{amsmath}
\usepackage{epsfig}
\usepackage{graphicx}

\usepackage{graphicx}
\usepackage{xypic}
\input xypic
\input xy
\xyoption{all}

\title{Rigidity of hyperbolic P-manifolds: a survey.}
\author{Jean-Fran\c{c}ois Lafont}
\address{Department of Mathematics,
The Ohio State University, 100 Math Tower, 231 West 18th Ave.,
Columbus, OH 43210} \email{jlafont@math.ohio-state.edu}

\newcommand{\mD}{\mathbb D}
\newcommand{\mH}{\mathbb H}

\newcommand{\G}{\Gamma}

\newcommand{\pinf}{\partial ^\infty}

\theoremstyle{proposition}
\newtheorem{Lem}{Lemma}[section]
\newtheorem{Prop}[Lem]{Proposition}

\newtheorem*{Def}{Definition}

\theoremstyle{plain}

\newtheorem{Thm}[Lem]{Theorem}

\newtheorem*{Thm2}{Main Theorem}

\theoremstyle{remark}

\begin{document}

\begin{abstract}
In this survey paper, we outline the proofs of the rigidity results
for simple, thick, hyperbolic P-manifolds found in \cite{L1},
\cite{L2}, \cite{L3}.
\end{abstract}

\maketitle

\section{Introduction.}

In recent years, there has been much interest in proving rigidity
type theorems for non-positively curved spaces.  All of these
results originated from Mostow's seminal work \cite{M}, in which he
showed that homotopy equivalent compact locally-symmetric spaces of
rank one are always isometric.

In the series of papers \cite{L1}, \cite{L2}, \cite{L3}, the author
exhibited rigidity results for a particularly nice class of singular
spaces.  The spaces in question, called hyperbolic P-manifolds, are
in some sense the simplest non-manifold CAT(-1) spaces one can
consider.  Let us start by defining these spaces:

\begin{Def}
A closed $n$-dimensional {\it piecewise manifold} (henceforth
abbreviated to P-manifold) is a topological space which has a
natural stratification into pieces which are manifolds. More
precisely, we define a $1$-dimensional P-manifold to be a finite
graph. An $n$-dimensional P-manifold ($n\geq 2$) is defined
inductively as a closed pair $X_{n-1}\subset X_n$ satisfying the
following conditions:

\begin{itemize}
\item Each connected component of $X_{n-1}$ is either an
$(n-1)$-dimensional P-manifold, or an $(n-1)$-dimensional manifold.
\item The closure of each connected component of $X_n-X_{n-1}$
is homeomorphic to a compact orientable $n$-manifold with boundary,
and the homeomorphism takes the component of $X_n-X_{n-1}$ to the
interior of the $n$-manifold with boundary; the closure of such a
component will be called a {\it chamber}.
\end{itemize}

\noindent Denoting the closures of the connected components of
$X_n-X_{n-1}$ by $W_i$, we observe that we have a natural map $\rho:
\coprod
\partial (cl(W_i))\longrightarrow X_{n-1}$ from the disjoint union
of the boundary components of the chambers to the subspace
$X_{n-1}$.  We also require this map to be surjective, and a
homeomorphism when restricted to each component. The P-manifold is
said to be \emph{thick} provided that each point in $X_{n-1}$ has at
least three pre-images under $\rho$. We will henceforth use a
superscript $X^n$ to refer to an $n$-dimensional P-manifold, and
will reserve the use of subscripts $X_{n-1},\ldots ,X_1$ to refer to
the lower dimensional strata.  For a thick $n$-dimensional
P-manifold, we will call the $X_{n-1}$ strata the {\it branching
locus} of the P-manifold.
\end{Def}

Intuitively, we can think of P-manifolds as being ``built'' by
gluing manifolds with boundary together along lower dimensional
pieces. Examples of P-manifolds include finite graphs and soap
bubble clusters.  Observe that compact manifolds can also be viewed
as (non-thick) P-manifolds. Less trivial examples can be constructed
more or less arbitrarily by finding families of manifolds with
homeomorphic boundary and glueing them together along the boundary
using arbitrary homeomorphisms. We now define the family of metrics
we are interested in.

\begin{Def}
A Riemannian metric on a 1-dimensional P-manifold (finite graph) is
merely a length function on the edge set.  A Riemannian metric on an
$n$-dimensional P-manifold $X^n$ is obtained by first building a
Riemannian metric on the $X_{n-1}$ subspace, then picking, for each
$cl(W_i)$ a Riemannian metric with totally geodesic boundary
satisfying that the gluing map $\rho$ is an isometry. We say that a
Riemannian metric on a P-manifold is hyperbolic if at each step, the
metric on each $cl(W_i)$ is hyperbolic.
\end{Def}

A hyperbolic P-manifold $X^n$ is automatically a locally $CAT(-1)$
space (see Chapter II.11 in Bridson-Haefliger \cite{BH}).
Furthermore, the lower dimensional strata $X_i$ are totally geodesic
subspaces of $X^n$.  In particular, the universal cover $\tilde X^n$
of a hyperbolic P-manifold $X^n$ is a $CAT(-1)$ space (so is
automatically $\delta$-hyperbolic), and has a well-defined boundary
at infinity $\partial ^\infty \tilde X^n$. Finally we note that the
fundamental group $\pi_1(X_n)$ is a $\delta$-hyperbolic group.  We
refer the reader to \cite{BH} for background on $CAT(-1)$ and
$\delta$-hyperbolic spaces.

\begin{Def}
We say that an $n$-dimensional P-manifold $X^n$ is {\it simple}
provided its codimension two strata is empty.  In other words, the
$(n-1)$-dimensional strata $X_{n-1}$ consists of a disjoint union of
$(n-1)$-dimensional manifolds.
\end{Def}

Examples of simple, thick, hyperbolic P-manifolds can easily be
constructed via arithmetic techniques.  The rigidity results
contained in the papers \cite{L1}, \cite{L2}, \cite{L3} can now be
summarized in the following:

\begin{Thm2}
Let $X^n$ be an $n$-dimensional, simple, thick, hyperbolic
P-manifold.  Then we have:
\begin{itemize}
\item if $n=2$, then $X^2$ is topologically rigid.
\item if $n\geq 3$, then $X^n$ is Mostow rigid.
\item if $n\geq 3$, then $\pi_1(X^n)$ is quasi-isometrically rigid.
\end{itemize}
In particular, if $n\geq 2$, the class of fundamental groups of
$n$-dimensional, simple, thick, hyperbolic P-manifolds exhibits
diagram rigidity.
\end{Thm2}

Let us briefly recall the terms involved.  By topological rigidity,
we mean that if $X^2_1$ and $X^2_2$ have isomorphic fundamental
group, then they are in fact homeomorphic.  By Mostow rigidity, we
mean that if $X^n_1$ and $X^n_2$ have isomorphic fundamental group,
then they are in fact isometric.  Quasi-isometric rigidity, refers
to the fact that the only groups which are quasi-isometric to
$\pi_1(X^n)$ are, up to finite extension, uniform lattices in
$Isom(\tilde X^n)$.  Finally, observe that the generalized
Seifert-Van Kampen theorem allows us to write the fundamental group
of any $n$-dimensional, simple, thick, hyperbolic P-manifold as the
direct limit of a canonical diagram of groups.  Diagram rigidity now
refers to the fact that two such direct limits are isomorphic if and
only if their diagrams are isomorphic.

The proof of these results naturally breaks down into two steps. The
first step requires an analysis of the topology of the boundary at
infinity of the universal cover of the simple, thick, hyperbolic
P-manifolds.  The second step `pushes back' the information garnered
at the boundary at infinity into the universal cover.

We will use the following standard notation: $S^n$ refers to an
$n$-dimensional sphere, $\mD^n$ to a closed $n$-dimensional ball,
and $\mD^n_\circ$ to an open $n$-dimensional ball.  $\mH^n$ refers
to the standard hyperbolic space of constant curvature $-1$.

\section{Topology of the boundary at infinity.}

We start by noting that, in a simple, thick, hyperbolic P-manifold
$X^n$, there is a codimension one singular set $B$ which consists of
a union of closed hyperbolic $(n-1)$-dimensional manifolds. Let
$\mathcal B$ denote the full pre-image of the singular set in the
universal cover $\tilde X^n$.  Observe that each connected component
$\tilde B_j$ of $\mathcal B$, lies as a totally geodesic subspace
(isometric to $\mH^{n-1}$) of the universal cover. This implies that
the boundary at infinity of $\tilde B_j$ (homeomorphic to $S^{n-2}$)
naturally embeds in the boundary at infinity of $\tilde X^n$.  By
abuse of notation, let $\pinf \mathcal B$ denote the subset of
$\pinf \tilde X^n$ consisting of the union of all the boundaries at
infinity of the various connected components of the branching locus,
i.e. $\pinf \mathcal B=\bigcup _j \pinf \tilde B_j$.

\begin{Thm}
Let $\pinf \tilde X^n_1$, $\pinf \tilde X^n_2$ ($n\geq 3$) be the
boundary at infinity of the universal cover of a pair of simple,
thick, hyperbolic P-manifold $X^n_i$, and let $\phi: \pinf \tilde
X^n_1 \rightarrow \pinf \tilde X^n_2$ be an arbitrary homeomorphism.
Then for each connected component $\tilde B_{1,j}$ in $\mathcal
B_1$, there is a connected component $\tilde B_{2,j^\prime}$ in
$\mathcal B_2$ with the property that $\phi(\pinf \tilde B_{1,j})=
\pinf \tilde B_{2,j^\prime}$.
\end{Thm}

In order to show this result, the natural starting point is to
obtain a topological characterization of points lying in $\pinf
\mathcal B$ that will distinguish them from points that are {\it
not} in $\pinf \mathcal B$.  This motivated the notion of
$n$-branching, defined as follows:

\begin{Def}
Define the 1-{\it tripod} $T$ to be the topological space obtained
by taking the join of a one point set with a three point set. Denote
by $*$ the point in $T$ corresponding to the one point set. We
define the $n$-{\it tripod} ($n\geq 2$) to be the space $T \times
\mD^{n-1}$, and call the subset $*\times \mathbb D^{n-1}$ the {\it
spine} of the tripod $T\times \mathbb D^{n-1}$.  The subset $*\times
\mathbb D^{n-1}$ separates $T\times \mathbb D^{n-1}$ into three open
sets, which we call the {\it open leaves} of the tripod.  The union
of an open leaf with the spine will be called a {\it closed leaf} of
the tripod. We say that a point $p$ in a topological space $X$ is
{\it $n$-branching} provided there is a topological embedding
$f:T\times \mathbb D^{n-1} \longrightarrow X$ such that $p\in
f(*\times \mathbb D^{n-1}_\circ)$.
\end{Def}

The topological characterization of points lying in $\mathcal B$ is
contained in the following:

\begin{Prop}
Let $X^n$ be a simple, thick, hyperbolic P-manifold ($n\geq 2$),
$\pinf \tilde X^n$ the boundary at infinity of its universal cover.
Then $p\in \pinf \tilde X^n$ is $(n-1)$-branching if and only if
$p\in \pinf \mathcal B$.
\end{Prop}

We will first explain how Theorem 2.1 follows from Proposition 2.2,
and then proceed to explain the argument behind Proposition 2.2.

\subsection{Reducing the Theorem to the Proposition}

Let us briefly explain how Theorem 2.1 can be obtained from
Proposition 2.2. We first observe that the property of being
$(n-1)$-branching is clearly a topological invariant.  So the
proposition implies that, if $\phi:\pinf \tilde X^n_1\rightarrow
\pinf \tilde X^n_2$ is an arbitrary homomorphism, we have that
$\phi(\pinf \mathcal B_1)=\pinf \mathcal B_2$.

Now recall that $\mathcal B_i= \bigcup _j \pinf \tilde B_{i,j}$. It
is easy to see that each $\pinf \tilde B_{i,j}$ is a closed subset
of $\pinf \tilde X^n_i$, that they are pairwise disjoint, and
furthermore if $n\geq 3$ that each $\pinf \tilde B_{i,j}$ is
connected (being homeomorphic to $S^{n-2}$). So if $n\geq 3$, we can
now apply a result of Sierpinski \cite{S} that states the following:
let $X$ be an arbitrary topological space, $\{C_i\}$ a countable
collection of disjoint path connected closed subsets in $X$.  Then
the path connected components of $\cup C_i$ are precisely the
individual $C_i$.  Applying this to the subset $\pinf \mathcal B_i$
in $\pinf \tilde X^n_i$, we see that the path-connected components
of $\pinf \mathcal B_i$ are precisely the individual $\pinf \tilde
B_{i,j}$. This forces the map $\phi$ to take each $\pinf \tilde
B_{1,j}$ homeomorphically onto some $\pinf \tilde B_{2,j^\prime}$,
which is the statement in Theorem 2.1.

In the case where $n=2$, one needs to be a bit more careful.  The
problem here is that when $n=2$, the branching locus in $X^2_i$ is
$1$-dimensional, i.e. consists of a finite union of closed
geodesics.  This implies that each $\pinf \tilde B_{i,j}$ in fact
consists of just a pair of points (i.e. an $S^0$), corresponding to
the two endpoints of a geodesic ray (a lift of the closed geodesic).
In particular, in this case, each subset $\pinf \tilde B_{i,j}$ is
{\it not} connected, and so Sierpinski's result is of no use.  In
order to get around this, we refine the information we have in the
$2$-dimensional setting by considering {\it separation properties}
of pairs of points in $\pinf \tilde X^2_i$.  Namely, one shows the
following:

\begin{Prop}
Let $\pinf \tilde X^2$ be the boundary at infinity of the universal
cover of a simple, thick, hyperbolic P-manifold of dimension $2$.
Then for a pair of points $x,y\in \pinf \tilde X^2$, the following
statements are equivalent:
\begin{itemize}
\item $\pinf \tilde X^2 -\{x,y\}$ has at least three path-connected
components.
\item $\pinf \tilde B_j =\{x,y\}$ for some connected component $B_j$
in $\mathcal B$.
\end{itemize}
\end{Prop}

The argument for this Proposition relies on the fact that $\pinf
\tilde X^2$ contains many isometrically embedded $\mH^2$ (this is
ensured by the thickness hypothesis).  Each of these yields an
embedded $S^1$ in the boundary at infinity, allowing us to construct
paths between points in the boundary at infinity by concatenating
paths traveling along the various $S^1$'s.  The details of the proof
can be found in \cite{L3}, and constitute the main difference
between the $2$-dimensional case and the higher dimensional cases.

\subsection{The argument for the Proposition.}

We now proceed to explain how to prove Proposition 2.2.  Recall that
the Proposition states that points in $\pinf \tilde X^n$ are
$(n-1)$-branching precisely if they lie in some $\pinf \tilde B_j$.

Well one direction of the implication is relatively easy: to show
that each point in $\pinf \tilde B_j$ is $(n-1)$-branching, we start
by observing that $\tilde B_j$, being a connected component of the
lift of the branching locus, is in fact an isometrically embedded
$\mH^{n-1}$ (this follows immediately from the simplicity
assumption), and hence that $\pinf \tilde B_j$ is homeomorphic to an
$S^{n-2}$.

Furthermore, it is easy to show that there exist, inside $\tilde
X^n$, three isometrically embedded ``half $\mH ^n$'s'', with the
property that their common boundary is $\tilde B_j$ (this follows
from the thickness hypothesis).  The net effect is that, on the
level of the boundary at infinity, one can find three disjoint
embedded $(n-1)$-dimensional closed disks $\mD^{n-1}$ with the
property that their boundaries all map homeomorphically to the
$\pinf \tilde B_j$ (an embedded $S^{n-2}$). But this immediately
gives you that the points in $\pinf \tilde B_j$ are
$(n-1)$-branching, completing one of the desired implications.

The reverse implication is considerably harder.  One wants to argue
that if $p\in \pinf \tilde X^n$ does {\it not} lie in one of the
$\pinf \tilde B_j$, then $p$ is {\it not} $(n-1)$-branching.  So
pick a point that does not lie in $\pinf \tilde B_j$.  With some
work, one can argue that it is suffices to assume that $p$ lies in
$\pinf \tilde W$ where $\tilde W$ is a connected lift of one of the
chambers.  We now would like to prove that the point $p$ is not
$(n-1)$-branching, i.e. to show that there are {\bf no} injective
maps $f:T\times \mD^{n-2}\rightarrow \pinf \tilde X^n$ satisfying
$p\in f(*\times \mD^{n-2})$.

The immediate difficulty is that one does not know much about the
topology of $\pinf \tilde X^n$.  But one of the nice properties of
CAT(-1) spaces is the existence of continuous maps from the boundary
at infinity to the {\it link} of any point inside the space.  Recall
that the link of a point essentially encodes the space of
``directions'' at that point.  The key feature we will need is that
the link of a point in the interior of $\tilde W$ is homeomorphic to
$S^{n-1}$ (since such a point has a manifold neighborhood).  If
$\rho$ denotes the projection map, we can now study the composite
map:
$$\rho \circ f:T\times \mD^{n-2}\longrightarrow \pinf \tilde
X^n\longrightarrow S^{n-1}$$ Note that it is immediate that points
in $S^{n-1}$ cannot be $(n-1)$-branching (for instance, they have
the wrong local homology), and hence the composite map $\rho \circ
f$ must fail to be injective.

The trick now boils down to ensuring that the composite map fails to
be injective at a point $q$ in the target {\it where the projection
map $\rho$ is injective}.  This will immediately imply that the
original map $f$ was not injective at the point $\rho^{-1}(q)$.  So
the next step is to understand the subset of $S^{n-1}$ where the map
$\rho$ is injective.  This is not too hard: an easy exercise in
elementary hyperbolic geometry gives you that the set of injective
points in $S^{n-1}$ consists of the complement of a countable dense
set of disjoint open metric balls in $S^{n-1}$. Let us denote this
set by $I$, and the subset consisting of the boundaries of the
various open metric balls by $\partial I\subset I$ (note that this
subset is {\it not} the boundary of the set $I$ in the topological
sense).

One also has that the point $\rho(p)$ in fact lies in $I-\partial
I$. Points inside $I-\partial I\subset S^{n-1}$ can be seen to have
a very special property: if an open set $U\subset S^{n-1}$ contains
one of these points in its closure, then $U$ must in fact contain a
point in $I$.  Since the set $I$ consists of those points where
$\rho$ is injective, we have reduced the proof of Proposition 2.2 to
the following:

\vskip 5pt

\noindent {\bf Claim:} Let $\rho\circ f: T\times
\mD^{n-2}\rightarrow S^{n-1}$ be as above.  Then there exists an
open set $U\subset S^{n-1}$ with the property that:
\begin{itemize}
\item for every point $y$ in $U$, $(\rho \circ f)^{-1}(y)$ consists
of at least two points.
\item the closure of $U$ contains $\rho (p)$.
\end{itemize}

\vskip 5pt

In order to show the claim, we first establish that an entire
neighborhood $V$ of $f^{-1}(p)$ inside the spine $*\times \mD^{n-2}$
maps to the set $I$.  Now consider the restriction of the composite
map $\rho \circ f$ to each of the three closed leaves of $T\times
\mD^{n-2}$.  Note that each of the closed leaves is homeomorphic to
$\mD^{n-1}$, with boundary homeomorphic to $S^{n-2}$.  Hence the
restriction of the composite to the boundary of each closed leaf
yields a continuous map from $S^{n-2}$ into $S^{n-1}$.  Furthermore,
the three maps $F_1,F_2,F_3$ obtained from the three closed leaves
have the property that they coincide on the spine, and in particular
on the open subset $V$ of the spine.

To conclude the argument for the claim, we would like to make
precise the following heuristic argument:
\begin{enumerate}
\item since the composite $\rho \circ f$ is injective on $V$, its
image is a topological embedding of $\mD^{n-2}$ inside $S^{n-1}$,
and hence should ``locally separate'' $S^{n-1}$ into two components.
\item the restriction of $\rho\circ f$ to each of the three closed
leaves must surject onto one of the two components, and hence one of
the two components must lie in the image of two distinct closed
leaves.
\end{enumerate}

In trying to make the above argument precise, one ends up involved
in the study of the maps $F_1,F_2,F_3$, which are maps from
$S^{n-2}$ into $S^{n-1}$ with the property that they are injective
on the small open set $V$.  The key property concerning these maps is the
following strong form of the Jordan separation theorem:

\begin{Thm}
Let $f:S^{n-2}\rightarrow S^{n-1}$ be a continuous map, and assume that
$f$ is injective on an open set $V\subset S^{n-2}$.  Then:
\begin{itemize}
\item $f(S^{n-2})$ separates $S^{n-1}$.
\item there are precisely two connected components $U_1,U_2$ in $S^{n-1}
-f(S^{n-2})$ having the property that there closure intersects $f(V)$.
\item if $F:\mD^{n-1}\rightarrow S^{n-1}$ is any extension of the map $f$,
then $F$ must surject onto either $U_1$ or $U_2$.
\end{itemize}
\end{Thm}

In the first paper \cite{L1} in the series, Theorem 2.4 was shown in
the case where $n=3$.  Note that in this case, one is looking at
maps from $S^1$ to $S^2$; the situation is then greatly simplified
by the Schoenflies Theorem, which was the key to establishing
Theorem 2.4 in the $n=3$ case.

When $n>3$, we know that the Schoenflies Theorem fails, and the
argument used in \cite{L1} has no chance of extending to higher
dimensions. So in \cite{L2}, a different argument was used to
establish Theorem 2.4. The argument relies on Alexander duality,
local Betti numbers, as well as the sophisticated codimension one
taming results of Bing \cite{B}, Ancel-Cannon \cite{AC}, and Ancel
\cite{A}.

\section{Concluding the proofs.}

As we mentioned earlier, Theorem 2.1 is the key result towards
proving the Main Theorem.  We now proceed to explain how to conclude
the proof of the Main Theorem.

\subsection{The common feature.}

First of all, recall that quasi-isometries between
$\delta$-hyperbolic spaces induce homeomorphisms between their
boundaries at infinity.  Note that an isomorphism between
fundamental groups yield quasi-isometries of the universal cover.

Now Theorem 2.1 essentially tells you that if $X^n_i$ are a pair of
simple, thick, hyperbolic P-manifolds, and if $\G_i=\pi_1(X^n_i)$
are the two fundamental groups, then any isomorphism
$\phi:\G_1\rightarrow \G_2$ naturally induces a bijection from the
connected components of $\mathcal B_1$ to the connected components
of $\mathcal B_2$.

The next step is to also obtain a bijection between the {\it lifts
of chambers} in the respective $\tilde X^n_i$.  In order to do this,
we exploit the separation properties of the $\pinf B_{i,j}\subset
\pinf \tilde X^n_i$.  Depending on the dimension, we prove that:
\begin{itemize}
\item if $n\geq 3$, each path-connected component of
$\pinf \tilde X^n - \pinf \mathcal B$ containing at least two points
corresponds canonically with a unique connected lift of a chamber in
$\tilde X^n$.
\item if $n=2$, we define an equivalence relation $\equiv$ on $\pinf
\tilde X^2-\pinf \mathcal B$ by defining $x\equiv y$ provided $\pinf
\tilde X^2-\{x,y\}$ has precisely two connected components, and show
that there is a canonical correspondance between equivalence classes
$\mathcal C$ in $\pinf \tilde X^2$ satisfying $|\mathcal C|>1$ on
the one hand, and connected lifts of chambers in $\tilde X^2$ on the
other.
\end{itemize}
Furthermore, in both cases above, if $W$ is the connected lift of a
chamber associated to the subset $Y\subset \pinf \tilde X^n-\pinf
\mathcal B$ (where $Y$ is as above, according to whether $n=2$ or
$n\geq 3$), then we also have that $\pinf W\subset \pinf \tilde X^n$
coincides with the closure of $Y$ in $\pinf \tilde X^n$.

The upshot is that the homeomorphism between the boundaries at
infinity take boundaries of lifts of chambers in $\pinf \tilde
X^n_1$ to boundaries of lifts of chambers in $\pinf \tilde X^n_2$.

\subsection{Mostow and topological rigidity.}

At this point, we have that the isomorphism $\phi: \G_1\rightarrow
\G_2$ induces bijections between:
\begin{itemize}
\item connected components of $\mathcal B_1$ and connected
components of $\mathcal B_2$, and
\item connected lifts of chambers in $\tilde X^n_1$ and connected
lifts of chambers in $\tilde X^n_2$.
\end{itemize}
Furthermore, the bijections are {\it compatible}, in the sense that
if a connected component of $\mathcal B_1$ is contained in the
connected lift of a chamber in $\tilde X^n_1$, then the same
statement holds for the objects bijectively associated to them.

Next we observe that the boundaries at infinity $\pinf \tilde X^n_i$
also come equipped with an action of $\G_i=\pi_1(X^n_i)$ by
homeomorphisms, and that the induced homeomorphism $\phi^\infty:
\pinf \tilde X^n_1\rightarrow \pinf \tilde X^n_2$ is
$(\G_1,\G_2)$-equivariant. It is easy to see that this implies that
the bijection between the connected components of the $\mathcal B_i$
actually {\it descends} to a bijection between the connected
components of the respective branching locis in the $X^n_i$.
Similarly, we obtain a bijection between the chambers of $X^n_1$ and
those of $X^n_2$.  Furthermore, this pair of bijections are again
compatible.

Now note that the chambers and connected components of the branching
locus are actually topological spaces.  We would like to ensure that
the bijections we've obtained actually preserve the topology of the
objects we are looking at.  This is achieved by establishing that
the fundamental group of the chambers (and branching locis) can be
detected from the boundary at infinity.  The technical statement is
the following:

\begin{Lem}
The stabilizer of a connected lift of a chamber under the
$\G$-action by deck transformations coincides with the stabilizer of
its boundary at infinity under the induced $\G$-action on $\pinf
\tilde X^n$.  The same statement holds for connected lifts of the
branching locus.
\end{Lem}

\noindent The proof of this Lemma is not very hard, and makes use of
the fact that the $\G$-action on $\pinf \tilde X^n$ exhibits
sink/source dynamics.  Now the $(\G_1,\G_2)$-equivariance of the
homeomorphism between the boundaries at infinity immediately gives
that the bijections between chambers (and between branching locis)
preserve the respective fundamental groups.

If $n\geq 3$, we can now apply Mostow rigidity for hyperbolic
manifolds with totally geodesic boundaries (see Frigerio \cite{Fr1})
to conclude that the chambers in bijective correspondance are
actually isometric. A little more work ensures that the isometries
between chambers glue together to give a global isometry from
$X^n_1$ to $X^n_2$, completing the argument for Mostow rigidity when
$n\geq 3$.

If $n=2$, we recall that an oriented surface with boundary is
topologically determined by its fundamental group and the number of
boundary components it has.  Since we know that chambers in $X^2_1$
correspond bijectively to chambers in $X^2_2$, and that
corresponding chambers have the same fundamental group and the same
number of boundary components, we can (being a little careful) glue
together the homeomorphisms between the various chambers to obtain a
global homeomorphism.  This gives us the desired topological
rigidity result when $n=2$.

Finally, let us say a few words about diagram rigidity.  Note that
we can associate to each P-manifold a diagram of groups obtained as
follows:
\begin{itemize}
\item corresponding to chambers and connected components of the
branching locus we associate vertices labelled with the fundamental
group of the respective object.
\item corresponding to each containment of a connected component of
the branching locus in a chamber, we associate a directed edge
between the corresponding vertices, labelled with the morphism
between fundamental groups induced by the inclusion.
\end{itemize}
The generalized Seifert-Van Kampen theorem now tells us that the
fundamental group of the P-manifold is in fact the direct limit of
the diagram of groups described above.

But we know that if the dimension is $n\geq 2$, any abstract
isomorphism between fundamental groups of a pair of $n$-dimensional
simple, thick, hyperbolic P-manifolds is in fact induced by a
homeomorphism (an isometry if $n\geq 3$) between the respective
P-manifolds.  So in particular, if we have an abstract isomorphism
between the direct limits of two diagrams as above, there is in fact
an isomorphism between the underlying graphs of the diagrams of
groups, having the property that corresponding vertices have
isomorphic fundamental groups.  Furthermore, the isomorphisms
between the vertex groups can be chosen to be compatible with the
edge morphisms (up to inner automorphisms), which is precisely the
statement of diagram rigidity.

\subsection{Quasi-isometric rigidity.}

For quasi-isometric rigidity, we will appeal to the following
well-known result (see for instance the survey by Farb \cite{Fa}):

\begin{Lem}
Let $X$ be a proper geodesic metric space, and assume that every
quasi-isometry from $X$ to itself is in fact a bounded distance (in
the sup norm) from an isometry. Furthermore, assume that a finitely
generated group $G$ is quasi-isometric to $X$.  Then there exists a
cocompact lattice $\Gamma\subset Isom(X)$, and a finite group $F$
which fit into a short exact sequence:
$$0\longrightarrow F\longrightarrow G\longrightarrow \Gamma
\longrightarrow 0$$
\end{Lem}
\vskip 5pt

\noindent In view of the Lemma, all we need to establish is that
every quasi-isometry of $\tilde X^n$ is bounded distance from an
isometry. In order to do this, we note that a quasi-isometry $f$ of
$\tilde X^n$ induces a self-homeomorphism of $\pinf \tilde X^n$.
From our knowledge of the topology of $\pinf \tilde X^n$, this
implies that for every lift of a chamber in $\tilde X^n$, the
quasi-isometry $f$ maps it to within finite distance of a (possibly
different) lift of a chamber.

By suitably perturbing the map $f$ by a bounded amount, we can
assume that the restriction of $f$ maps each lift of a chamber
quasi-isometrically into the lift of another chamber.  At this
point, we appeal to a well known `folk-theorem' (a rigorous proof of
which can be found in \cite{Fr2}): if $Y_1$,$Y_2$ are the universal
covers of two compact hyperbolic manifolds with non-empty totally
geodesic boundaries, and $g:Y_1\rightarrow Y_2$ is a quasi-isometry,
then there is an isometry $\bar g:Y_1\rightarrow Y_2$ within bounded
distance of $g$.

One now patches together the isometries between the various chambers
to obtain that the original $f:\tilde X^n\rightarrow \tilde X^n$ is
at finite distance from an isometry.  There are two points to be
careful with:
\begin{itemize}
\item we need to ensure that the isometries on the chambers do
indeed glue together to give a global isometry.
\item we need to make sure that the resulting isometry is at bounded
distance from the map $f$ (i.e. that the bounded distance on each of
the chambers is actually uniformly bounded).
\end{itemize}
But neither of these two points are very hard to establish.
Applying the Fact now completes the argument for quasi-isometric
rigidity.

\section{Concluding remarks.}

We conclude this paper by pointing out to the reader that there are
several natural questions which are still unanswered:

\begin{itemize}
\item Does topological rigidity hold for simple, thick, {\it negatively} curved
P-manifolds of dimension $\geq 5$?
\item For which classes of diagrams of groups does diagram rigidity
hold?
\item Does Mostow type rigidity still hold if we remove the
simplicity hypothesis (i.e. allow for a more complicated singular
set)?
\end{itemize}

For a more thorough discussion of these questions, we refer the
reader to the problem list in \cite{L4}.


\begin{thebibliography}{99}

\bibitem[A]{A} F.D. Ancel, `Resolving wild embeddings of codimension-one manifolds in
manifolds of dimensions greater than $3$', Special volume in honor
of R. H. Bing (1914--1986), {\it Topology Appl.} 24 (1986), pp.
13--40.

\bibitem[AC]{AC} F.D. Ancel and  J.W. Cannon,
`The locally flat approximation of cell-like embedding relations',
{\it Ann. of Math.} (2) 109 (1979), pp. 61--86.

\bibitem[B]{B} R.H. Bing,
`Approximating surfaces by polyhedral ones', {\it Ann. of Math.} (2)
65 (1957), pp. 465--483.

\bibitem[BH]{BH} M.R. Bridson and A. Haefliger,
{\it Metric spaces of non-positive curvature} (Springer-Verlag,
Berlin, 1999).

\bibitem[Fa]{Fa} B. Farb, `The quasi-isometry classification of
lattices in semisimple Lie groups, {\it Math. Res. Lett.} 4 (1997),
pp. 705-717.

\bibitem[Fr1]{Fr1} Frigerio R., `Hyperbolic manifolds with geodesic boundary which are determined
by their fundamental group', {\it Topology Appl.} 145 (2004), pp.
69--81.

\bibitem[Fr2]{Fr2} Frigerio R., `Commensurability of hyperbolic manifolds with geodesic boundary', preprint
available on the ArXiv at
http://front.math.ucdavis.edu/math.GT/0502209.

\bibitem[L1]{L1} J.-F. Lafont, `Rigidity result for certain 3-dimensional singular spaces
and their fundamental groups', {\it Geom. Dedicata} 109 (2004), pp.
197--219.

\bibitem[L2]{L2} J.-F. Lafont, `Strong Jordan separation and applications to rigidity',
to appear in {\it J. London Math. Soc.}

\bibitem[L3]{L3} J.-F. Lafont, `Rigidity of geometric amalgamations of free
groups', preprint available on the ArXiv at
http://front.math.ucdavis.edu/math.GR/0506518

\bibitem[L4]{L4} J.-F. Lafont, `Some open problems in Geometry and
Topology', in preparation.

\bibitem[M]{M} G.D. Mostow, {\it Strong rigidity of locally symmetric spaces} (Princeton University Press,
Princeton, N.J., 1973).

\bibitem[S]{S} W. Sierpinski, `Un th\'eor\`eme sur les continus', {\it Tohoku Math. Journ.} 13 (1918),
pp. 300-303.

\end{thebibliography}
\end{document}